\documentstyle[11pt,mm,amsmath,amssymb,graphics,url,epsfig]{article}

\newtheorem{theorem}{Theorem}[section]

\newtheorem{algorithm}[theorem]{Algorithm}
\newtheorem{remark}[theorem]{Remark}

\def\citet{\cite}
\newcommand{\comment}[1]{}

\newcommand{\osty}{\displaystyle}
\newcommand{\xreal}{\xi}

\newcommand{\xvar}{X}
\newcommand{\zvar}{z}

\newcommand{\QQ}{{\mathbb Q}}
\newcommand{\RR}{{\mathbb R}}

\def\startpage{66}  
\def\endpage{73}   
\def\vol{28}       
\def\year{2009}    
\def\runauthor{Feng Guo}

\def\runtitle{SDPTools: High Precision SDP Solver in Maple}

\begin{document}

\title{ SDPTools: High Precision SDP Solver in Maple}

\author{Feng Guo\\
 Key Laboratory of Mathematics Mechanization\\
 Institute of Systems Science,   AMSS,
 Academia Sinica\\
 Beijing 100190,  China\\
 fguo@mmrc.iss.ac.cn}
 \date{}

\maketitle

\begin{abstract}
Semidefinite programs are an important class of convex optimization problems.
It can be solved efficiently by SDP solvers in Matlab, such as SeDuMi, SDPT3,
DSDP. However, since we are running fixed precision SDP solvers in Matlab, for
some applications, due to the numerical error, we can not get good results.
SDPTools is a Maple package to solve SDP in high precision. We apply SDPTools
to the certification of the global optimum of rational functions. For the
Rump’s Model Problem, we obtain the best numerical results so far.
\end{abstract}


\section{Preliminaries}
\label{sec-0}

SDPTools is a Maple package to solve  a semidefinite program (SDP):
\begin{equation}\label{sdpprimal}
\begin{aligned}
\mathrm{minimize} & \quad c^Tx   \\
~\mbox{s.t.}~ & \quad  \mathit{F}(x) \ge 0,
\end{aligned}
\end{equation}
where
\begin{equation*}
\mathit{F}(x) \overset{\triangle}= \mathit{F_0} + \sum^{m}_{i=1}
x_i\mathit{F_i}.
\end{equation*}
The problem data are the vector $c \in {\RR}^m$ and $m+1$ symmetric
matrices $F_0, \dots, F_m \in {\RR}^{n \times n}$. The dual problem
associated with the semidefinite program~(\ref{sdpprimal}) is
\begin{equation}\label{sdpdual}
\begin{aligned}
\mathrm{maximize} &  \quad \mathrm{-Tr} \mathit{F_0}Z \\
~\mbox{s.t.}~ &  \quad \mathrm{Tr} \mathit{F_i}Z=c_i, \quad i=1,...,m, \\
 & \quad Z \ge 0.
\end{aligned}
\end{equation}
Here the variable is the matrix $Z=Z^T \in {\RR}^{n \times n}$,
which is subject to $m$ equality constraints and the matrix
nonnegativity condition.

The SDP~(\ref{sdpprimal}) and its dual problem~(\ref{sdpdual}) can
be solved efficiently by algorithms  in SeDuMi~\cite{Sturm99},
SDPT3~\cite{Toh98sdpt3}, DSDP~\cite{dsdp5}, and
SDPA~\cite{Fujisawa98sdpa}. However, since we are running fixed
precision SDP solvers in Matlab, we can only obtain  numerical
positive semidefinite matrices  which satisfy equality or inequality
constrains approximately. For some applications, such as Rump's
model problem~\cite{Rump06}, due to the numerical error, the
computed lower bounds can even be significantly larger than upper
bounds, see, e.g., Table~1 in~\cite{KLYZ08}. These numerical
problems motivate us to consider how to use symbolic computation
tools such as Maple to obtain SDP solutions with high accuracy.

\section{An algorithm  for solving SDPs}
\label{sec-2}

Semidefinite programs are an important class of convex optimization
problem for which readily computable self-concordant barrier
functions are known, so interior-point methods are applicable. Our
algorithm for solving SDPs is based on the potential reduction
methods mentioned in \cite{vandenberghe96semidefinite}, which is to
minimize the potential function (55) in
\cite{vandenberghe96semidefinite}
\begin{eqnarray}
\varphi(x,Z)&\overset{\triangle}= &  v \sqrt{n} \text{log}({\mathbf Tr}F(x)Z)+\psi(x,Z) \\
                              & = & (n+v \sqrt{n})\text{log}({\mathbf Tr}F(x)Z)-\text{log det} F(x)-\text{log det}Z -n\text{log} n \nonumber
\end{eqnarray}
The first term in the right hand of the first equation is the
duality gap of the pair $x, Z$ and $\psi(x,Z)$ denotes the deviation
from the centrality. When we do iterations to minimize
$\varphi(x,Z)$ with a strict feasible start point, the first term
insures $(x_k, Z_k)$ approach the optimal point and the second
guarantees all $(x_k, Z_k)$ are feasible. More details, see
\cite{vandenberghe96semidefinite}. A general outline of a potential
reduction method is as follows

\begin{algorithm}\label{PotentialReduction}PotentialReduction

\noindent
 \textbf{Input:} Strictly feasible
$x$, $Z$, and a tolerance $\epsilon$.

\noindent
 \textbf{Output:} Updated strictly feasible
$x$ and $Z$.

\noindent
 \textbf{Repeat:}
\begin{enumerate}

\item \label{step:1}  Find suitable
directions $\delta x$ and $\delta Z$.

\item \label{step:2}  Find $p, q\in
\RR$ that minimize $\varphi(x+p\delta x, Z+q\delta Z)$.

\item  \label{step:3} Update: $x:=x+p \delta x$ and $Z:=Z+q \delta Z$.

\end{enumerate}
\textbf{until} duality gap $\le \epsilon$.
\end{algorithm}

 In  Step~\ref{step:1}, an
obvious way to compute search directions $\delta x$ and $\delta Z$
is to apply Newton's method to $\varphi$. However the potential
function $\varphi$ is not a convex function, since the first term:
$(n+v \sqrt{n} )\text{log}({\mathbf Tr}F_0Z + c^Tx) $ is concave in
$x$ and $Z$ and hence contributes a negative semidefinite term to
Hessian of $\varphi$. We adapt potential reduction method 2
mentioned in \cite{vandenberghe96semidefinite}, which is based on
the primal system only. Namely, we choose direction $\delta x$ that
minimize a quadratic approximation of $\varphi (x+v,Z)$ over all $v
\in \RR^m $. In order to apply the Newton's method, the second
derivative of the concave term is ignored. It is equivalent to solve
the following linear equations:
\begin{equation}\label{direc-x}
\begin{aligned}
F\delta Z^p F+ \sum_{i=1}^m \delta x_i^p F_i = - \rho FZF +F \\
{\mathbf Tr}F_j \delta Z^p =0, \quad j=1,\cdots, m.
\end{aligned}
\end{equation}
we choose $\delta Z^p$ as the dual search direction, see
\cite{vandenberghe96semidefinite}. Practically, this method seems
perform better than method 1 which treats the primal and dual
semidefinite program symmetrically.

In Step~\ref{step:2}, we use plane search to choose the lengths of
the steps made in the directions $\delta x$ and $\delta Z$, see
\cite{vandenberghe96semidefinite}.
\begin{equation}\label{length}
\begin{aligned}
\text{minimize}\quad  &(n+ v \sqrt{n})\text{log}(1+c_1p+c_2q)- \sum_{i=1}^n\text{log}(1+p \mu_i)- \sum_{i=1}^n\text{log}(1+q\nu_i) \\
~{\mbox s.t.}~ \quad &p_\text{min} \le p \le p_\text{max}, \quad
q_\text{min} \le q \le q_\text{max}
\end{aligned}
\end{equation}
where
\[ c_1=\frac{c^T\delta x}{{\mathbf Tr}F(x)Z}, \quad   c_2=\frac{{\mathbf Tr}F_0\delta Z}{{\mathbf
Tr}F(x)Z},
\]
$\mu_1, \dots, \mu_n$ and $\nu_1, \dots, \nu_n$ are the
eigenvalues of $F^{-1/2}\delta FF^{-1/2}$ and $Z^{-1/2}\delta
ZZ^{-1/2}$, respectively.

Since the objective is a quasiconvex function, we combine the
Newton's direction and the steepest decent direction together to
minimize it. If the Hessian matrix at an iteration point is positive
definite, we use the former, otherwise we use the latter. After we
get the decent direction $\delta p$ and $\delta q$, we use bisection
method to compute the one dimensional search of $r$, and update
$p_{k+1}=p_k+r \cdot \delta p, \ q_{k+1}=q_k+r \cdot \delta q$,
accordingly.

\begin{figure}
\begin{center}
\scalebox{0.5}{\includegraphics{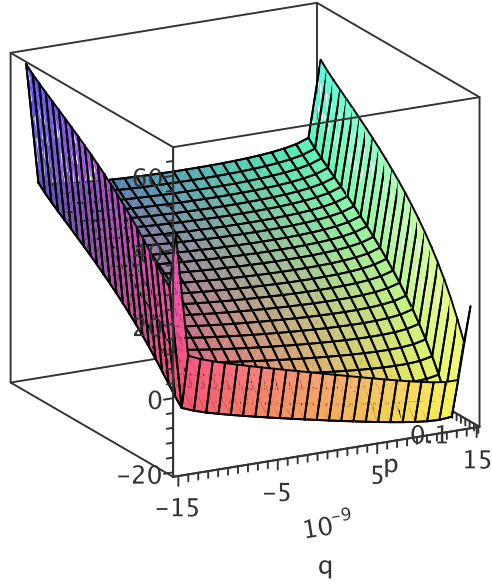}}
\caption{\label{figure_1}The objective function for plane search}
\end{center}
\end{figure}

\begin{remark}
 Since the objective function in (\ref{length}) is very sensitive at the
optimizer(see Figure \ref{figure_1}), the operation must be
performed carefully. A tiny numerical error can cause endless loops.
In order to get higher accuracy result by doing more iterations, we
need to set larger digits.
\end{remark}

\begin{remark}
We should also pay attention to choose the initial point $(p, q)$.
Usually, we choose $(0,0)$ and it works well. But there comes
troubles for Rump's problem(see, \ref{sec-rump}). When $n=8$, the
iteration for one dimensional search does not converge. Because our
methods only converge locally, if $(0,0)$ is too far from the
optimizer, the iterations may not converge. If we denote the optimal
solution as $(p^*,q^*)$, according to our experiments, we find $p^*$
is always attained near the largest feasible value and $q^*$ near
$0$ (see Figure \ref{figure_1}). It is always reasonable since when
computing search direction $\delta x$ and $\delta Z$, we choose
method 2, which is based on the primal system only. So we only
compute the approximate Newton's direction of $x$ and choose
accessary result $\delta Z^p$ of (\ref{direc-x}) as the dual search
direction. So the step length $q$ for $\delta Z$ should be near $0$.
So we choose $(\widehat{p}/2, 0)$ as the initial point in plane
search where $\widehat{p}$ denotes the largest feasible value of
$p$, and it solves the trouble of Rump's problem when $n=8$.
\end{remark}

The algorithm {\sf PotentialReduction} needs initial strictly
feasible primal and dual points. We adapt the Big-M method when
neither a strictly feasible $x$ nor a strictly feasible $Z$ is
known. After modification, the primal problem becomes
\begin{equation}\label{primal}
\begin{aligned}
\mathrm{minimize} & \quad c^Tx+\mathit{M_1}t  \\
s.t & \quad  \mathit{F}(x)+t\mathit{I} \ge 0, \\
&\quad \mathbf{Tr}\mathit{F} \le \mathit{M_2}, \\
& \quad t \ge 0
\end{aligned}
\end{equation}
The dual becomes
\begin{equation}\label{dual}
\begin{aligned}
\mathrm{maximize} &\quad \mathbf{-Tr} \mathit{F_0}(Z-z\mathit{I})-\mathit{M_2}z \\
s.t & \quad\mathbf{Tr} \mathit{F_i}(Z-z\mathit{I})=c_i, \quad i=1,...,m. \\
&\quad\mathbf{Tr}\mathit{Z} \le \mathit{M_1}, \\
&\quad \mathit{Z} \ge 0,z \ge 0
\end{aligned}
\end{equation}
If $M_1$ and $M_2$ are big enough, this entails no loss of
generality(assuming the primal and dual problem are bounded). After
modification, problem (\ref{primal}) and (\ref{dual}) can be written
as (\ref{sdpprimal}) and (\ref{sdpdual}). Then it is easy to compute
its strictly feasible points, see \cite{vandenberghe96semidefinite}

A brief description of the main functions contained in SDPTools is
following:

\begin{itemize}

\item
{\sffamily Solution\_{}of\_{}Directions} \ solves the directions
$\delta x$ and $\delta Z$ from (\ref{direc-x}).

\item
{\sffamily Expression\_{}of\_{}PlaneSearch} \ gets the objective of
(\ref{length}).

\item
{\sffamily StepsLength\_{}PlaneSearch} \ solves problem
(\ref{length}) and get step length $p, q$ for $x, Z$ respectively.

\item
{\sffamily BigM\_{}Case3} \ transforms (\ref{primal}), (\ref{dual})
into (\ref{sdpprimal}), (\ref{sdpdual}) respectively, and computes a
strict feasible start point $(x_0, Z_0)$.

\item
{\sffamily Solve\_{}SDP\_{}Method2} \
solves the SDP
(\ref{sdpprimal}) and (\ref{sdpdual}).

\end{itemize}

\section{Certified Global Optimum of Rational Functions}
\label{sec-2}

In SDPTools, we apply the SDP algorithm (\ref{PotentialReduction})
to compute and certify the lower bounds of rational functions with
rational coefficients.
\begin{align}\label{rationalfun}
\min_{\xreal\in \RR^n} \frac{f(\xreal)}{g(\xreal)}\text{\quad (where
$g(\xreal) > 0$ for all $\xreal\in\RR^n$)}
\end{align}
where $f(\xreal), g(\xreal) \in \QQ^n$.  The number   $r$ is a lower
bound of  (\ref{rationalfun}) if and only if the polynomial $f(x)-r
g(x)$ is nonnegative.  Therefore we focus on the following
minimization of a rational function by SOS:
\begin{equation}\label{SOS}\left.
\begin{aligned}
r^* := & \sup_{r \in \RR,W} r \\
s.t & \ f(x)-r g(x) = m_{d}(x)^T \cdot W \cdot m_{d}(x) \\
& W \succeq 0, \quad W^T=W
\end{aligned} \right\}
\end{equation}
where $m_d(x)$ is the column vector of all terms in $X_1, \dots,
X_n$ up to degree $d$. A detailed description of SOS relaxations and
its dual problems, see \cite{Nie}. The problem (\ref{SOS}) is a
semidefinite program and can be written as (\ref{sdpdual}). With
packages for solving SDP, we can obtain a numerical positive
semidefinite matrix $W$ and floating point number $r^*$ which
satisfy approximately:
\begin{equation}
f(x)-r^* g(x) \approx  m_{d}(x)^T \cdot W \cdot m_{d}(x),\
W\succapprox 0
\end{equation}
To certify $r^*$, we convert $r^*, W$ to rational ones and project
orthogonally $W$ onto the affine linear hyperplane:
\begin{equation}
\begin{aligned}
\chi = \{ A| A^T &= A, f(x)-r g(x)= m_{d}(x)^T \cdot A \cdot
m_{d}(x),\quad \text{for some $r$} \}
\end{aligned}
\end{equation}
and hope $W$ to be positive semidefinte after projection.

However, if we run the fixed precision SDP solver in Matlab, the
output $W$ are too coarse to be projected into the cone of positive
semidefinite matrices. In \cite{KLYZ08} they refined the $r^*, W$ by
Gauss-Newton method. Here our high precision SDP solver  shows its
advantage. It is implemented in Maple, which has arbitrarily high
precision, so it can compute $r^*, W$ with high accuracy. Without
refinement, the projection can be done successfully.

\begin{figure}
\begin{center}
\scalebox{0.5}{\includegraphics{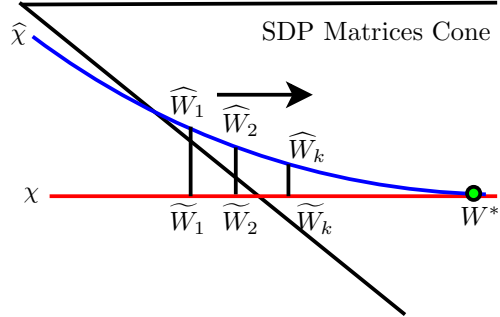}}
\put(-90,105){\small{SDP Matrices Cone}} \put(-15,35){\small{$W^*$}}
\put(-77,33){\small{$\widetilde{W}_k$}}
\put(-105,33){\small{$\widetilde{W}_2$}}
\put(-125,33){\small{$\widetilde{W}_1$}}
\put(-125,77){\small{$\widehat{W}_1$}}
\put(-105,70){\small{$\widehat{W}_2$}}
\put(-80,60){\small{$\widehat{W}_k$}}
\put(-180,45){\footnotesize{$\chi$}}
\put(-185,105){\footnotesize{$\widehat{\chi}$}}
\caption{\label{figure_2}Rationalization of SOS}
\end{center}
\end{figure}

Usually, it is not easy to find a strictly feasible point for
(\ref{SOS}) and we need the Big-M method to avoid it. After convert
(\ref{SOS}) to the form (\ref{dual}), the SOS relaxation of
(\ref{rationalfun}) becomes
\begin{equation}\label{MSOS}\left.
\begin{aligned}
\widehat{r}^* := \sup_{\widehat{r} \in \RR,\widehat{W}} \quad & \widehat{r} -M_2 z \\
s.t. \quad & f(x)-\widehat{r} g(x)+z(m_{d}(x)^T \cdot m_{d}(x)) = m_{d}(x)^T \cdot \widehat{W} \cdot m_{d}(x), \\
& \widehat{W} \succeq 0, \quad \widehat{W}^T=\widehat{W}, \quad z
\ge 0
\end{aligned} \right\}
\end{equation}
It is obvious to see that problem (\ref{SOS}) and (\ref{MSOS}) are
equivalent. And the variable $z$ is required to be $0$ at the
optimizer of problem(\ref{MSOS}). We define
\begin{equation}
\begin{aligned}
\widehat{\chi} = \{ \widehat{A}\ |\  \widehat{A}^T &= \widehat{A}, f(x)-\widehat{r} g(x)+z(m_{d}(x)^T \cdot m_{d}(x)) \\
&= m_{d}(x)^T \cdot \widehat{A} \cdot m_{d}(x),\quad \text{for some
$\widehat{r}, z$} \}
\end{aligned}
\end{equation}
Note that $\widehat{\chi}$ will meet $\chi$ at the optimizer
$(W^*,r^*)$, see Figure \ref{figure_2}.

After each iteration, we get a point $(\widehat{r}, \widehat{W})$ on
$\widehat{\chi}$. Then we convert $(\widehat{r}, \widehat{W})$ to
rational ones and project matrix $\widehat{W}$ onto $\chi$,  and
denote the nearest matrix in $\chi$ by $\widetilde{W}$. Because the
problem(\ref{SOS}) and (\ref{MSOS}) are equivalent, we can expect
that after the first few ones, at each iteration we can obtain a
positive semidefinite matrix $\widetilde{W}$ and a certified $r$.
Different from the method in \cite{KLYZ08} which only certifies a
given $r$, we can get series of certified lower bounds $r_n$. The
more iterations we do, the better $r_n$ we get. The process is shown
in Figure \ref{figure_2}.

In order to reduce the problem size, we also exploit sparsity. Given
a polynomial $p(x)=\sum_{\alpha}p_{\alpha}x^{\alpha}$, the cage of
$p$, is the convex hull of $\sup(p)= \{\alpha |\ p_{\alpha} \neq 0
\}$.
\begin{theorem}(\cite{Reznick78})
For a polynomial $p$, $C(p^2)=2C(p)$; for any positive semidefinite
polynomial $f$ and $g$, $C(f) \subseteq C(f+g)$; if $f=\sum_jg_j^2$
then $C(g_j) \subseteq \frac{1}{2}C(f)$.
\end{theorem}
With the above theorem, we can remove the redundant monomials in
(\ref{SOS}).

SDPTools also has functions to compute the convex hull of given
finite points in $n$ dimensional vector space. They are mainly based
on the Quickhull algorithm in \cite{quickhull} which assumes the
affine dimension of the input points is also $n$. If the affine
dimension, denoted by $r$, is less than $n$, SDPTools provides a
function to   compute $r$ and affinely transforms the input points
to a $r$ dimensional vector space where we can use the Quickhull
algorithm. For a given point,  SDPTools has a function to judge
whether it is in the convex hull.

SDPTools  provides a function to  compute the monomials appearing in
(\ref{SOS}). For example, when solving Rump's problem, the monomials
appearing in $m_d(x)$ are obtained by proof in \cite{KLYZ08} while
we  can find them automatically by SDPTools!

The followings are the main functions to compute and certify the
lower bound of a given rational function:

\begin{itemize}
\item
{\sffamily affineDims} \ computes the affine dimension of input
points.
\item
{\sffamily affineTrans} \ affinely transforms the given points to a
$r$ dimensional vector space.
\item
{\sffamily convexHull} \ computes the convex hull of given points.
\item
{\sffamily inConvexHull} \ judges whether $p$ is in the convex hull
defined by the points in $S$.
\item
{\sffamily getSDP} \ transforms problem (\ref{SOS}) as form
(\ref{sdpdual}).
\item
{\sffamily projSOS} \ converts $r, W$ to rational ones and projects
$W$ onto $\chi$.
\item
{\sffamily certifiSOS} \ computes and certifies the lower bound of
an input rational function.
\end{itemize}

\section{A numerical experiment: solving Rump's Model Problem}
\label{sec-rump}

The following introduction to Rump's model problem (see,
\cite{Rump06}) is mostly from \cite{KLYZ08}. This problem is related
to structured condition numbers of Toeplitz matrices and polynomial
factor coefficient bounds, asks for $n=1,2,3,\ldots $ to compute the
global minima
\begin{equation*}
\begin{array}{r@{}@{}l@{}@{}l}
\mu_n = {} &\osty \min_{P,\,Q}\ &\osty \frac{\|PQ\|_2^2}{\|P\|_2^2\|Q\|_2^2} 
\\[2ex]
           &\text{s.\ t.\ } &\osty P(\zvar)= \sum_{i=1}^{n} p_i \zvar^{i-1} 
                            , Q(\zvar)=\sum_{i=1}^{n} q_i \zvar^{i-1} \in \RR[\zvar] \setminus \{0\}.
\end{array}
\end{equation*}
It has been shown in~\cite{RumpSek06} that polynomials $P,Q$
realizing the polynomials achieving $\mu_n$ must be symmetric
(self-reciprocal) or skew-symmetric. Thus the problem can be
rewritten into three optimization problems with three different
constraints
\begin{equation*}
\begin{array}{rlll}
k=1\colon &p_{n+1-i}=p_i, & q_{n+1-i}=q_i, &1 \leq i \leq n, \\
k=2\colon &p_{n+1-i}=p_i, & q_{n+1-i}=-q_i,&1 \leq i \leq n, \\
k=3\colon & p_{n+1-i}=-p_i,& q_{n+1-i}=-q_i,&1 \leq i \leq n,
\end{array}
\end{equation*}
and the smallest of three minima is equal to $\mu_n$. For all three
cases, we minimize the rational function $f(\xvar)/g(\xvar)$ with
\begin{equation*}
\begin{array}{l@{}@{}l}
  f(\xvar) &
\displaystyle {} = \|P Q\|_2^2=\sum_{k=2}^{2n} (\sum_{i+j=k} p_i
q_j)^2,
\quad
 g(\xvar) 
\displaystyle {} = \|P\|_2^2 \|Q\|_2^2=(\sum_{i=1}^{n}
p_i^2)(\sum_{j=1}^{n} q_j^2)
\end{array}
\end{equation*}
and the variables
$ \xvar=\{p_1,\ldots,p_{n(P)}\}\cup\{q_1,\ldots,q_{n(Q)}\}, $
where  $n(P)=n(Q)=\lceil {n}/{2} \rceil$.

In this paper, we use higher precision SDP solver in SDPTools  to
solve Rump's model problem and obtain much better certified lower
bounds, see Table \ref{table}.

\begin{table}
\centerline{\footnotesize
\begin{tabular}{|c|c|c|c|r|l|l|} \hline
$n$ & $k$ & \# iter & prec. &
\multicolumn{1}{|c|}{secs/iter}&\multicolumn{1}{|c|}{lower bound
$r_n$}& \multicolumn{1}{c}{upper bound}
\\ \hline 4 & 2 & 50 &
4$\times$ 15& 0.71   & 0.01742917332143265287    &
0.01742917332143265289 \\\hline 5 & 1 & 50 & 4 $\times$ 15& 2.03   &
0.00233959554815559112    & 0.00233959554815559113 \\ \hline 6 & 2 &
50 & 4 $\times$ 15& 1.76   & 0.00028973187527968191    &
0.00028973187527968193  \\ \hline 7 & 1 & 75 & 5 $\times$ 15& 11.36
& 0.00003418506980008284    & 0.00003418506980008285   \\ \hline 8 &
2 & 75 & 5 $\times$ 15& 12.49  & 0.00000390543564975572    &
0.00000390543564975573    \\ \hline 9 & 1 & 75 & 5 $\times$ 15&
84.12  & 0.43600165391810484612e-06& 0.43600165391810484613e-06 \\
\hline 10& 2 & 75 & 5 $\times$ 15& 92.79  &
0.47839395687709759326e-07& 0.47839395687709759327e-07 \\ \hline 11&
1 & 85 & 5 $\times$ 15& 622.03 & 0.51787490974469905330e-08&
0.51787490974469908331e-08   \\ \hline 12& 2 & 85 & 5 $\times$ 15&
634.48 & 0.55458818311631347611e-09&0.55458818311631347612e-09\\
\hline 13& 1 & 100& 5 $\times$ 15& 3800.0 &
0.58866880811866093129e-10 & 0.58866880811866093130e-10 \\ \hline
14& 2 & 100& 5 $\times$ 15& 3800.0 & 0.62024449920539050219e-11
&0.62024449920539050220e-11 \\ \hline 15& 1 & 120& 6 $\times$ 15&
15000  & 0.64943654185809512879e-12 &0.64943654185809512880e-12  \\
\hline 16& 2 & 120& 6 $\times$ 15& 23000  &
0.67636042558221379057e-13 &0.67636042558221379058e-13 \\ \hline
\end{tabular}
} \caption{\label{table}The certified lower bounds}
\end{table}

\section{Conclusion}
SDPTools is a Maple package having the following functions: to solve
a general SDP, to compute and certify the lower bounds of rational
functions. It is also a tool to compute the convex hull of given
finite points in order to explore the sparsity structure.

SDPTools is still under development, we hope to implement more
efficient algorithms to solve SDP, to give upper bounds for $M_1,
M_2$ when applying Big-M method in SOS relaxation (\ref{MSOS}), to
detect the infeasibility for SDP and to explore more sparsity
structures for problems having large size.

\def\refname{\Large\bfseries References}
\bibliographystyle{plain}
\bibliography{fguo}

\end{document}